\documentclass[]{article}  
\usepackage{CJK}
\usepackage{amssymb}
\usepackage{bbm}
\usepackage{amsfonts}
\usepackage{mathrsfs}
\usepackage{graphicx}
\usepackage{amsmath}
\usepackage{amsthm}
\usepackage{xypic}
\usepackage{graphicx}
\usepackage{times}
\usepackage{geometry}
\usepackage{color}
\usepackage{indentfirst}
\usepackage{cases}
\usepackage{hyperref}
\usepackage{graphicx}
\usepackage{float}
\usepackage{amsthm}
\usepackage{amssymb}
\usepackage{amsmath}
\usepackage{mathrsfs}
\usepackage{indentfirst}
\usepackage{geometry}
\usepackage{authblk}
\usepackage{hyperref}

\renewcommand\figurename{Fig}
\newtheorem{Theorem}{Theorem}[section]
\theoremstyle{definition}

\newtheorem{Lemma}[Theorem]{Lemma}
\newtheorem{Proposition}[Theorem]{Proposition}

\theoremstyle{definition}

\title{A Characterization for a graph an eigenvalue of multiplicity $2c(G) + q_s(G) - 1$}
\author{Songnian Xu, Wenhao Zhen, Dein Wong\thanks{Corresponding author. E-mail address: wongdein@163.com.}}
\affil{\textit{School of Mathematics, China University of Mining and Tecnology, Xuzhou, China.}}
\date{}

\begin{document}
\baselineskip 17pt

\title{A Characterization for a graph an eigenvalue of multiplicity $2c(G) + q_s(G) - 1$}

\author{Songnian Xu\\
{\small  Department of Mathematics, China University of Mining and Technology}\\
{\small Xuzhou, 221116, P.R. China}\\
{\small E-mail: xsn1318191@cumt.edu.cn}\\ \\
Dein Wong\thanks{Corresponding author}\\
{\small Department of Mathematics, China University of Mining and Technology}\\
{\small Xuzhou 221116, P.R. China}\\
{\small E-mail: }}

\date{}
\maketitle

\begin{abstract}
Let $G$ be a simple connected graph and  $P = u_1 u_2 \ldots u_t u_{t+1}$ be an induced subgraph of $G$.
If $d_G(u_1) = 1$, $d_G(u_2) = \ldots = d_G(u_t) = 2$, and $d_G(u_{t+1}) \geq 3$, then we refer to $u_1 u_2 \ldots u_t$ as a pendant path $P_t$.
For $G \in \mathbb{G}_s=\{G\mid$ every pendent path in $G$ is at least of order $s\}$, let $Q_s(G)$ be the set of vertices in $G$ that are distance $s$ from some pendant vertices, and let $|Q_s(G)| = q_s(G)$. 
For $G \in \mathbb{G}_s$, Li et al. (2024) proved that when $\lambda$ is not an eigenvalue of $P_s$ and $G$ is neither a cycle nor a starlike tree $T_k$, it holds that $m_G(\lambda) \leq 2c(G) + q_s(G) - 1$ and characterized the extremal graphs when $G$ is a tree.
In this article, we characterize the extremal graphs for which  $m_G(\lambda) = 2c(G) + q_s(G) - 1$ when $G \in \mathbb{G}_{s}$ and $\lambda\notin \sigma(P_s)$.
\end{abstract}

\let\thefootnoteorig\thefootnote
\renewcommand{\thefootnote}{\empty}
\footnotetext{Keywords: eigenvalue multiplicity; cyclomatic number; extremal graphs}

\section{Introduction}
In this paper, we restrict our attention to simple, connected, and finite graphs.
A simple undirected graph $G$ is denoted by $G = (V(G), E(G))$, where $V(G)$ represents the vertex set and $E(G)$ represents the edge set.
We denote that $G$ is a tree by $G = T$.
A graph $H$ is defined as a subgraph of the graph $G$ if $V(H) \subseteq V(G)$ and $E(H) \subseteq E(G)$.
Furthermore, $H$ is termed an induced subgraph of $G$ if any two vertices in $V(H)$ are adjacent in $H$ if and only if they are adjacent in $G $; this relation is denoted as $H \leq G$.
The order of $G$, representing the number of vertices, is denoted by $|G|$.
For a vertex $x \in V(G)$ and a subgraph $H \leq G$, let $N_H(x) = \{ u \in V(H) \mid uv \in E(G) \}$.
If $K \subseteq V(G)$, we denote the subgraph of $G$ induced by $K$ as $G[K]$.
We sometimes express $G - K$ or $G - G[K]$ to signify $G[V(G) \setminus V(K)]$.
A vertex $v$ in $G$ is classified as a pendant vertex if $d_G(v) = 1$, where $d_G(v)$ denotes the number of adjacent vertices of $v$ in $V(G)$.
If $d_G(u) \geq 3$, we denote $u$ as a high degree vertex in $G$.
Let $M_G$ represent the set of all high degree vertices in $G$.
For a connected graph $G$ with at least one cycle, a maximal leaf-free induced subgraph of $G$ is called the $plinth$ of $G$, which is written as $G^{\circ}$.
Let $\sigma(G)$ denote the set of all eigenvalues of the adjacency matrix $A(G)$ of the graph $G$ and $f_G(\lambda)=|\lambda I-A(G)|$.

Let $P = u_1 u_2 \ldots u_t u_{t+1}$ be an induced subgraph of $G$.
If $d_G(u_1) = 1$, $d_G(u_2) = \ldots = d_G(u_t) = 2$, and $d_G(u_{t+1}) \geq 3$, then we refer to $u_1 u_2 \ldots u_t$ as a pendant path $P_t$.
For $s\leq t+1$, the vertex $u_s$ is called a $(s-1)$-pendant vertex.
Let $Q_s(G)$ denote the set of all $s$-pendant vertices of $G$, and define $q_s(G) = |Q_s(G)|$.

In 2020, Wang et al. \cite{WL1} proved that when $G$ is not a cycle, there holds $m_G(\lambda) \leq 2c(G) + p(G) - 1$, where $\lambda \in \mathbb{R}$.
In 2024, Zhang et al. \cite{ZZ1} provided a complete characterization of graphs for which $m_G(\lambda) = 2c(G) + p(G) - 1$.
Subsequently, let $\parallel\lambda\parallel$ is the absolute absolute value of $\lambda$, Zhen et al. \cite{ZD1} generalized the aforementioned result by demonstrating that when $\parallel\lambda\parallel \geq 2$, it holds that $m_G(\lambda) \leq c(G) + q(G)$, and they characterized the extremal graphs when $G$ is a cactus or a leaf-free graph with at most one vertex of degree 3.

Let $\mathbb{G}_s=\{G\mid$ every pendent path in $G$ is at least of order $s\}$.
In 2024, for $G \in \mathbb{G}_s$ and $\lambda\notin \sigma(P_s)$, Li et al. \cite{LW1} proved that $m_G(\lambda) \leq 2c(G) + q_s(G)$ and established that equality holds if and only if $G$ is a cycle or a starlike tree $T_k$.
Consequently, when $G$ is neither a cycle nor a starlike tree $T_k$, it follows that $m_G(\lambda) \leq 2c(G) + q_s(G) - 1$; Li et al. also characterized the extremal graphs when $G$ is a tree and equality holds.
In this article, we provide a complete characterization of the graphs for which $m_G(\lambda) = 2c(G) + q_s(G) - 1$ when $G \in \mathbb{G}_s$ and $\lambda\notin \sigma(P_s)$.

\section{Our main result}

Before presenting the main results, we will first provide definitions of some concepts that will be helpful for the subsequent writing of this paper.

$T_k$: The tree $T_k\in \mathbb{G}_s$ on $ks + 1$ vertices is obtained from $k$ copies of path $P_s$ by connecting a pendant vertex of each path with a new vertex.
In particular, the tree $T_2$ is the path $P_{2s+1}$.

$T_{k_1,k_2,l}$: For positive integers $k_1$, $k_2$ and $l$, denote by $T_{k_1,k_2,l}$ the tree obtained from two trees $T_{k_1}$ and $T_{k_2}$ by identifying their centers with two distinct pendant vertices of a path of order $l\geq 2$, respectively.
Note that $T_{k_1,k_2,l}\in \mathbb{G}_s$.

$T_{k,l}$: When $l \geq 1$ is an integer, let $T_{k,l}$ be the tree on $sk + l$ vertices obtained from $T_k$ by identifying its center with a
pendant vertex of $P_l$.
We call another pendant vertex of $P_l$ the $l$-vertex of $T_{k,l}$.
In particular, when $l = 0$, we define $T_{k,0}$ as the union of $k$ disjoint paths $P_s$.

$C_{k,l}$: When $l \geq 1$ is an integer, let $C_{k,l}$ be the graph on $k + l-1$ vertices obtained from $C_k$ by identifying its a vertex with a
pendant vertex of $P_l$.
We call another pendant vertex of $P_l$ the $l$-vertex of $C_{k,l}$.

$s$-$p$-$deletion$: Let $H$, $T_k$, and $P_l$ be three pairwise disjoint graphs, where $u \in V(H)$ and $d_H(u) \geq 2$, and $v$ is the central vertex of $T_k$.
Identify $u$ and $v$ with two pendant vertices of $P_l$ respectively (and when $l = 1$, identifying $u$ and $v$), we obtain the graph $G=H_u+P_l(u,v)+(T_k)_v$.
Then, we say that $H$ is obtained from $G$ through an $s$-$p$-deletion process.

$c$-$p$-$deletion$: Let $H$, $C_k$, and $P_l$ be three pairwise disjoint graphs, where $u \in V(H)$ and $d_H(u) \geq 2$, and $v$ is a vertex of $C_k$.
Identify $u$ and $v$ with two pendant vertices of $P_l$ respectively (and when $l = 1$, identifying $u$ and $v$), we obtain the graph $G=H_u+P_l(u,v)+(C_k)_v$.
Then, we say that $H$ is obtained from $G$ through an $c$-$p$-deletion process.

Before presenting our main results, we first provide some key results from \cite{LW1} and \cite{ZZ1}, which will be beneficial for describing and proving our findings.

\begin{Theorem}\cite{ZZ1}
Let $G$ be a connected graph and $\lambda$ be a real number.

Then $m_G(\lambda)=2c(G)+p(G)-1$ if and only if $G$ is one of the following forms:

($i$) $\lambda=2$; $G$ is a cycle.

($ii$) $\lambda=-2$; $G$ is a cycle with even order.

($iii$) $G$ is the bicyclic graph $B(l,1,k)$, where $l,k$ are both multiple of $4$.

($iv$) $G$ is the bicyclic graph $\theta(l',x',k')$, where $l'=x'=k'\equiv 0$ or $2 \ (\text{mod} \ 4)$.

($v$) $\lambda=2cos\frac{i\pi}{m+1}$, where $i$, $m+1$ are co-prime integers with $1\leq i\leq m$; $G$ is a tree that satisfies the following conditions: $d(u,v)\equiv m \ (\text{mod} \ m+1)$ for any $u\in P_G$ and $v\in M_G$ (if exists); $d(u_1,u_2)+1\equiv m \ (\bmod \ m+1)$
 for any two distinct pendant vertices $u_1$, $u_2$ in $G$.

 ($vi$) $\lambda=2cos\frac{i\pi}{m+1}$, where $i$, $m+1$ are co-prime integers with $1\leq i\leq m$; $G$ is obtained from a tree $T$ with $m_T(\lambda)=p(T)-1$ by turning $c(G)\geq1$ pendant $P_m$ into $c(G)$ pendant cycles, each cycle has $\lambda$ as its eigenvalue.
\end{Theorem}

\begin{Theorem}\cite[2.12]{LW1}
Let $G\in \mathbb{G}_s$ be a connected graph on $n$ vertices with cyclomatic number $c(G)$.
For any eigenvalue $\lambda$ of $G$, if $\lambda\notin \sigma(P_s)$, then
$$m_G(\lambda)\leq 2c(G)+q_s(G).$$
The equality holds if and only if $G$ is either a starlike tree $T_k$ on $ks+1$ vertices with $\lambda\in \sigma(T_k)\setminus \sigma(P_s)$, or a cycle $C_n$ with  $\lambda\in \sigma(C_n)\setminus (\sigma(P_s)\cup\{\pm2\})$.
\end{Theorem}

\begin{Theorem}\cite[3.1]{LW1}
Let $T\in \mathbb{G}_s\setminus T_{\frac{n-1}{s}}$ be a tree of order $n$.
$\lambda\notin \sigma(P_s)$ is an eigenvalue of $T$, then

$$m_T(\lambda)\leq q_s(T)-1$$

The equality holds if and only if $T$ is one of the following graphs:

(1) the tree $T_{k_1, k_2,l}$ with $\lambda\in \sigma(T)\setminus \sigma(P_s)$ and $n = (k_1 + k_2)s + l$;

(2)a tree $T_0$ on $n_0$ vertices with $q_s\geq 3$ pendant vertices $u_1,\ldots,u_{q_s}$, by identifying the center of $T_{k_i}$ with the pendant vertex $u_i$ for each $i\in \{1,\ldots,q_s\}$, where $n=n_0+s\Sigma^{q_s}_{i=1}k_i$.
Meanwhile, for $i\in \{1,\ldots, q_s\}$, we have $\lambda$ is an eigenvalue of pendant $T_{k_i,l_i}$, any pair of vertices in $V(T_0)$ of degree at least 3 are non-adjacent, and all internal paths of $T_0$ share the eigenvalue $\lambda$.
\end{Theorem}

Let $H$ be a graph, $u \in V(H)$ such that $d_H(u) \geq 2$, and $v$ is an $l$-vertex of $T_{k,l}$.
If $G$ is the graph obtained by adding an edge between $u$ and $v$, then we refer to $T_{k,l}$ as a $pendant$-$T_{k,l}$ of $G$.
Similarly, if $v$ is an $l$-vertex of $C_{g,l}$ and $G$ is the graph obtained by adding an edge between $u$ and $v$, then we refer to $C_{g,l}$ as a $pendant$-$C_{g,l}$ of $G$.

For $G\in \mathbb{G}_s$ and $\lambda \notin \sigma(P_s)$, when $G$ is a tree or a leaf-free graph, Theorems 2.1 and 2.3 have already addressed all cases where $m_G(\lambda) = 2c(G) + q_s(G) - 1$.
Therefore, in this article, we do not need to consider the situation where $G$ is leaf-free graph or a tree. Below, we present main results of this paper.

The following Theorems 2.4 describe the situation when $c(G) = q_s(G) = 1$.

\begin{Theorem}
Let $G=C_u+P_l(u,v)+(T_k)_v$,$\lambda\notin \sigma(P_s)$ and $l\geq1$.
Then $m_G(\lambda) = 2c(G) + q_s(G) - 1=2$ if and only if $G$ is one of the following forms:

(1) $l=1$,  $k=1$, $m_C(\lambda)=2$ and $m_{P_{s-1}}(\lambda)=1$.

(2) $l\geq3$, $m_C(\lambda)=2$ and $m_{T_{k,l-2}}(\lambda)=1$.
\end{Theorem}

Let $m$ be a positive integer, it's well known that:

(1) The eigenvalues of $P_m$ are $\{2cos\frac{i\pi}{m+1}|i=1,2,\ldots,m\}$;

This imply that, if there are two co-prime integers $i$, $m+1$ with $1\leq i\leq m$ such that $\lambda=2cos\frac{i\pi}{m+1}$, $\lambda$ is an eigenvalue of $P_k$ if and only if $k\equiv m \ (\bmod \ m+1))$.

(2) The eigenvalues of $C_m$ are $\{2cos\frac{2i\pi}{m}|i=0,1,2,\ldots,m-1\}$.

Let $T \in \mathbb{G}_s\setminus T_{\frac{n-1}{s}}$, $\lambda \in \sigma(T)\backslash \sigma(P_s)$, and $m_T(\lambda) = q_s - 1$.
Replace $h \leq q_s$ pendant-$T_{k_i,l_i}$ in $T$ with pendant-$C_{g_i,f_i}$ to obtain the graph $G$ for $1\leq i\leq m$, where $m_{C_{g_i,f_i}}(\lambda) = 2$ (and according to Theorem 2.1, this is equivalent to letting $\lambda = 2 \cos \frac{i \pi}{m+1}$, where $i$ and $m + 1$ are coprime for $1 \leq i \leq m$, $\lambda \in \sigma(C_{g_i})$, and $f_i \equiv 1 \ (\bmod \ m+1))$.
The set of all such graphs $G$ is denoted by $\mathbb{B}(T,h)$.

The following Theorem 2.5 describes the situation when $c(G) + q_s(G) \geq 3$ and $m_G(\lambda) = 2c(G) + q_s(G) - 1$.
\begin{Theorem}
Let $G \in \mathbb{G}_{s}$, $G \neq G^{\circ}$ or $T$, $c(G)+q_s(G)\geq3$ and $\lambda \notin \sigma(P_s)$.
Then $m_G(\lambda) = 2c(G) + q_s(G) - 1$ if and only if $G \in \mathbb{B}(T, c(G))$.

\end{Theorem}

The following Proposition 2.6 indicates that the sufficiency of Theorem 2.4 and 2.5 holds.
Thus, in the subsequent sections, it suffices to prove the necessity.

\begin{Proposition}
Let $G$ is a graph that satisfies condition of Theorem 2.4 or 2.5, then $m_G(\lambda) = 2c(G) + q_s(G) - 1$.
\end{Proposition}

\begin{proof}

\textbf{Case 1}: When $G$ satisfies condition of Theorem 2.4, let $u_1 \in N_G(u) \setminus V(C)$.
According to the Interlacing Theorem, we have
$$m_G(\lambda) \geq m_{G-u_1}(\lambda) - 1 = m_C(\lambda) +1 - 1 = m_C(\lambda) = 2.$$
Furthermore, since $m_G(\lambda) \leq 2c(G) + q_s(G) - 1 = 2$ by Theorem 2.2, it follows that $m_G(\lambda) = 2$.

\textbf{Case 2}: When $G$ satisfies condition of Theorem 2.5, let $T_0$ be the tree corresponding to $T$ in Theorem 2.3(2).
According to the interlacing theorem, we have $m_G(\lambda) \geq m_{G - M_{T_0}}(\lambda) - |M_{T_0}|$.
Here, $G - M_{T_0}$ is the disjoint union of all pendant-$C_{g_i, l_i}$ and pendant-$T_{k_j, l_j}$ from $G$, along with all internal paths $P_{t_f}$ of $T_0$ (noting that there are $|M_{T_0}| - 1$ internal paths in total).
From the definition of $\mathbb{B}(T, c(G))$ and the results of Theorem 2.3(2), we know that $m_{C_{g_i, h_i}}(\lambda) = 2$, $m_{T_{k_i, l_i}}(\lambda) = 1$, and $m_{P_{t_f}}(\lambda) = 1$.
Therefore, we can conclude that
$$m_G(\lambda) \geq m_{G - M_{T_0}}(\lambda) - |M_{T_0}| = 2c(G) + q_s(G) + |M_{T_0}| - 1 - |M_{T_0}| = 2c(G) + q_s(G) - 1.$$

On the other hand, according to Theorem 2.2, we also know that
$$m_G(\lambda) \leq 2c(G) + q_s(G) - 1.$$
Thus, the proof is complete.
\end{proof}

\section{Proof of main result}

All graphs in this paper are simple undirected graphs. The adjacency matrix $A(G)$ of $G$ is an $n\times n$ square matrix whose $(i, j)$ entry takes 1 if vertices $i$ and $j$ are adjacent in $G$, and it takes 0 otherwise.
The eigenvalues of $A(G)$ are directly called the eigenvalues of $G$. For $\lambda\in \mathbb{R}$, let
 $$\mathbb{V}^{\lambda}_G=\{\alpha\in \mathbb{R}_n| A(G)\alpha=\lambda \alpha\}$$
 $$\mathbb{Z}^{X}_G= \{\alpha\in \mathbb{R}_n| \alpha_x = 0, x\in X\}$$

\begin{Lemma}\cite[2.2]{ZD1}
 Let $Y\subseteq X$ be subsets of $V(G)$ with $|X|-|Y| = m$.
 Then dim$(\mathbb{V}^{\lambda}_G\cap \mathbb{Z}^{Y}_G)\leq$ dim$(\mathbb{V}^{\lambda}_G\cap \mathbb{Z}^{X}_G) + m$.
 Furthermore, $m_G(\lambda)\leq$  dim$(\mathbb{V}^{\lambda}_G\cap \mathbb{Z}^{X}_G) +|X|$.
\end{Lemma}

\begin{Lemma}\cite[2.5]{LW1}
 For integers $k\geq1$ and $l\geq2$, if $\lambda\in \sigma(T_{k,l})$ and $\lambda\notin \sigma(P_s)$, then $\lambda\notin \sigma(T_{k,l-1})$.
 Moreover, $m_{T_{k,l}}(\lambda)=1$.
\end{Lemma}

Since $\lambda \notin \sigma(P_s)$, it is evident that the aforementioned Lemma 3.2 holds for the case when $l = 1$ by Lemma 3.1.
Therefore, in the subsequent sections of this paper, we will only require $l \geq 1$ when referring to Lemma 3.2.

\begin{Lemma}
Let $G$ be a connected graph, $\lambda \in \sigma(G)$ and $\lambda \notin \sigma(P_s)$. 
If $G$ is one of the following two cases, then we have $m_G(\lambda) \leq 2$.

(1) $G=C_u+P_l(u,v)+(T_k)_v$ and $l\geq1$;

(2) $G=C_{g,l}$ and $l\geq 1$.
\end{Lemma}
\begin{proof}
We will only prove case (1), as case (2) can be demonstrated in a similar manner.
Let $C$ be the unique cycle on $G$, and let $u \in V(C) \cap M_G$ and $z$ be a neighbor of $u$ on $C$.
According to Lemma 2.1, we have
$$m_G(\lambda) \leq \dim(\mathbb{V}^{\lambda}_G \cap \mathbb{Z}^{\{u,z\}}_G) + 2.$$
For $(x_1, x_2, \ldots, x_n)^{T} = \alpha \in \mathbb{V}^{\lambda}_G \cap \mathbb{Z}^{\{u,z\}}_G$, considering any pendant $P_s$ of $G$, we have $A\alpha = \lambda\alpha$ and  $\lambda x_i = \sum_{j \sim i} x_j$, which implies $\alpha\mid_{P_s} \in \mathbb{V}^{\lambda}_{P_s}$.
However, since we know that $\lambda \notin \sigma(P_s)$, it follows that $\alpha\mid_{P_s} = 0$.
Therefore, based on $A\alpha = \lambda\alpha$ and $\lambda x_i = \sum_{j \sim i} x_j$, we can conclude that $\alpha = 0$, which leads to $ m_{G}(\lambda)\leq2$.
\end{proof}

A vertex $u$ of the graph $G$ is call $D$-vertex (or ``downer") in $G$ for an eigenvalue $\lambda$ of $A(G)$ if $m_{A(G)}(\lambda)=m_{A(G-u)}(\lambda)+1$.
\begin{Lemma}\cite[2.3]{TK1}
Let $G$ be a connected graph, $A=A(G)$, and $\lambda\in \sigma(A)$.
Then, vertex $i$ is downer for $\lambda$ if and only if column (resp. row) $i$ of $A-\lambda I$ is a linear combination of the other columns (resp. rows) of $A-\lambda I$.
\end{Lemma}

The following lemma is essential for the proof of our theorem.
\begin{Lemma}
 Let $GuvH$ be the graph obtained from $G\cup H$ by adding an edge joining the vertex $u$ of $G$ to the vertex $v$ of $H$.
If $m_G(\lambda)=m_{G-u}(\lambda)+1$, then $m_{GuvH}(\lambda)=m_{G-u}(\lambda)+m_{H-v}(\lambda)$.
In particular, when $\lambda \in \sigma(G)$ but $\lambda \notin \sigma(G - u)$, we have $m_{GuvH}(\lambda) = m_{H-v}(\lambda)$.

\end{Lemma}

\begin{proof}
Without loss of generality, we consider the case when $\lambda = 0$.
$A(G)=\bordermatrix{
& G-u & u \cr
 &A(G-u) & \alpha  \cr
 &\alpha^{T}&  0
}
.$
According to Lemma 3.4, there exists an invertible matrix $P_1$ such that
$P_1A(G)P^{T}_1=\bordermatrix{
& G-u & u \cr
 &A(G-u) & 0  \cr
 &0&  0
}
$.
Let $P=
\begin{pmatrix}
  P_1& 0 \\
 0&  E_{|H|}
\end{pmatrix}
$,
Then we have $$PA(GuvH)P^{T}
  =\bordermatrix{
& G-u & u & v& H-v \cr
 &A(G-u) & 0 & 0 &0  \cr
 &0 & 0 & 1 &0 \cr
 &0 & 1 & 0 & \beta^{T} \cr
 &0 & 0  & \beta& A(H-v)
}.$$

Let $Q=\begin{pmatrix}
  E_{|G|-1}& 0&0&0 \\
 0&  1& 0& 0 \\
 0&0&1&0\\
 0&-\beta&0&E_{|H|-1}
\end{pmatrix}
$, then we have $QPA(GuvH)P^{T}Q^{T}=\begin{pmatrix}
  A(G-u)& 0&0&0 \\
 0&  0& 1& 0 \\
 0&1&0&0\\
 0&0&0&A(H-v)
\end{pmatrix}$.

Therefore, we have $m_{GuvH}(0) = n - \text{rank}(A(G - u)) - 2 - \text{rank}(H - v)$, which implies that $m_{GuvH}(0) = m_{G - u}(0) + m_{H - v}(0)$.
\end{proof}

\begin{Lemma}\cite{SD1}
Let $u$ be a pendant vertex and $v$ the quasi-pendant vertex adjacent to $u$ of $G$.
Then
$$f_G(\lambda)=\lambda f_{G-u}(\lambda)-f_{G-u-v}(\lambda).$$
\end{Lemma}

\textbf{Proof of Theorem 2.4}
\begin{proof}
When $G=C_u+P_l(u,v)+(T_k)_v$, $l\geq1$ and $m_G(\lambda)=2$.
By Lemma 3.1, we have $m_G(\lambda) \leq \dim(\mathbb{V}^{\lambda}_G \cap \mathbb{Z}^{v}_G) + 1$.
This implies that $\dim(\mathbb{V}^{\lambda}_G \cap \mathbb{Z}^{v}_G)\geq1$.
For any $(x_1, x_2, \ldots, x_n)^{T} = \alpha \in \mathbb{V}^{\lambda}_G \cap \mathbb{Z}^{v}_G$, based on $A\alpha = \lambda\alpha$ and $\lambda x_i = \sum_{i\sim j}x_j$, we find that $\alpha\mid_{P_s} \in \mathbb{V}^{\lambda}_{P_s}$ for any pendant path $P_s$ of $T_k$.
However, according to the problem statement, $\lambda \notin \sigma(P_s)$, which leads to $\alpha\mid_{P_s}  = 0$.
Thus, $\alpha\mid_C \in \mathbb{V}^{\lambda}_C$, and since $x_u = 0$ by $\alpha\mid_{P_s}  = 0$ and $\lambda x_i = \sum_{i\sim j}x_j$, we obtain from $A\alpha = \lambda\alpha$ that $\alpha\mid_{P_{|C|-1}} \in \mathbb{V}^{\lambda}_{P_{|C|-1}}$.
Therefore, we have $\lambda \in \sigma(P_{|C|-1}) \cap \sigma(C)$, which implies $m_C(\lambda) = 2$ and $m_{P_{|C|-1}}(\lambda) = 1$.

when $l=1$, we will prove that $k=1$.
If $k = 2$, let $P_s = v_1 v_2 \ldots v_s$ be any pendant path on $T_k$, where $d_G(v_1) = 1$, $d_G(v_2) = \ldots = d_G(v_s) = 2$, and $v_s \sim u$. 
Since $m_G(\lambda) = 2$, we know that $m_{G - v_1}(\lambda) \geq 1$, which implies $f_{G - v_1}(\lambda) = 0$. 
According to Lemma 3.6, we have $f_{G - v_1 - v_2}(\lambda) = \lambda f_{G - v_1}(\lambda) - f_G(\lambda) = 0$. 
Continuing this reasoning, we obtain $f_{G - v_1 - v_2 - \ldots - v_s - u}(\lambda) = 0$. 
Since $f_{G - v_1 - v_2 - \ldots - v_s - u}(\lambda) = f_{C - u}(\lambda) + (k - 1)f_{P_s}(\lambda)$, and knowing that $\lambda \in \sigma(C - u)$, we infer $f_{P_s}(\lambda) = 0$, which contradicts $\lambda \notin \sigma(P_s)$. 
Therefore, we conclude that $k = 1$ when $l=1$.
When $k = 1$, since $u$ is a $D$-vertex of $C$ with respect to $\lambda$, by Lemma 3.5 we have $m_G(\lambda) = m_{C - u}(\lambda) + m_{P_{s-1}}(\lambda) = 2$. 
Consequently, given $m_{C - u}(\lambda), m_{P_{s-1}}(\lambda) \leq 1$, we conclude that $m_{C - u}(\lambda) = m_{P_{s-1}}(\lambda) = 1$.

When $l = 2$, since $u$ is a $D$-vertex of $C$ with respect to $\lambda$, by Lemma 3.5, we have $m_G(\lambda) = m_{C-u}(\lambda) + k m_{P_s}(\lambda) = m_{C-u}(\lambda) = 1$.
This is in contradiction with $m_G(\lambda) = 2$, thus $l \neq 2$.

When $l \geq 3$, similarly, we have $m_G(\lambda) = m_{C-u}(\lambda) + m_{T_{k,l-2}}(\lambda) = 2$, which implies $m_{T_{k,l-2}}(\lambda) = 1$. Theorem 2.4 is proven.

\end{proof}

For $G \in \mathbb{G}_s$, if $m_G(\lambda) = 2c(G) + q_s(G) - 1$, we say that the graph $G$ is $\lambda$-$optimal$.

\begin{Lemma}
Let $G \in \mathbb{G}_s$, $\lambda\notin \sigma(P_s)$ and $q_s(G) \geq 1$.
If $G'$ is a graph obtained from $G$ through an $s$-$p$-deletion process and $G' \neq C$ or $T_k$, then if $G$ is $\lambda$-optimal, $G'$ is also $\lambda$-optimal.
\end{Lemma}

\begin{proof}
Without loss of generality, let $G = G'_u + P_l(u,v) + (T_k)_v$.
By Lemma 2.1, we know that $m_G(\lambda) \leq \dim(\mathbb{V}^{\lambda}_G \cap \mathbb{Z}^{v}_G) + 1$.
For any $(x_1,x_2,\ldots,x_n)^{T}=\alpha \in \mathbb{V}^{\lambda}_G \cap \mathbb{Z}^{v}_G$, from $A\alpha = \lambda \alpha$, we can deduce that for any pendant path $P_s$ of $T_k$, we have $\alpha\mid_{P_s} \in \mathbb{V}^{\lambda}_{P_s}$ .
However, since $\lambda \notin \sigma(P_s)$, it follows that $\alpha\mid_{P_s} = 0$.
Thus, from $\lambda x_i = \sum_{i\sim j}x_j$, we have $x_i = 0$ for $i\notin V(G')$, which implies that $\alpha\mid_{G'} \in \mathbb{V}^{\lambda}_{G'}$.
Consequently, we obtain
$$\dim(\mathbb{V}^{\lambda}_G \cap \mathbb{Z}^{v}_G) \leq m_{G'}(\lambda).$$
This leads to the conclusion that
$$m_{G'}(\lambda) \geq m_G(\lambda)-1 = 2c(G) + q_s(G) - 2.$$
Furthermore, since $c(G) = c(G')$ and $q_s(G) = q_s(G')+1$, it follows that $m_{G'} \geq 2c(G') + q_s(G') - 1$.
Combining these results and taking into account that $G' \neq C$ or $T_k$, we find that $m_{G'} \leq 2c(G') + q_s(G') - 1$ by Theorem 2.2.
Therefore, $G'$ is $\lambda$-optimal.
\end{proof}

Similar to Lemma 3.7, we present Lemma 3.8. 
Since The proof of Lemma 3.8 is similar to that of Lemma 3.7., we shall omit the proof of Lemma 3.8.

\begin{Lemma}
Let $G \in \mathbb{G}_s$, $\lambda\notin \sigma(P_s)$ and $q_s(G) \geq 1$.
If $G'$ is a graph obtained from $G$ through an $c$-$p$-deletion process and $G' \neq C$ or $T_k$, then if $G$ is $\lambda$-optimal, $G'$ is also $\lambda$-optimal.
\end{Lemma}

\begin{figure}[H]
  \centering
  \includegraphics[width=1.0\linewidth]{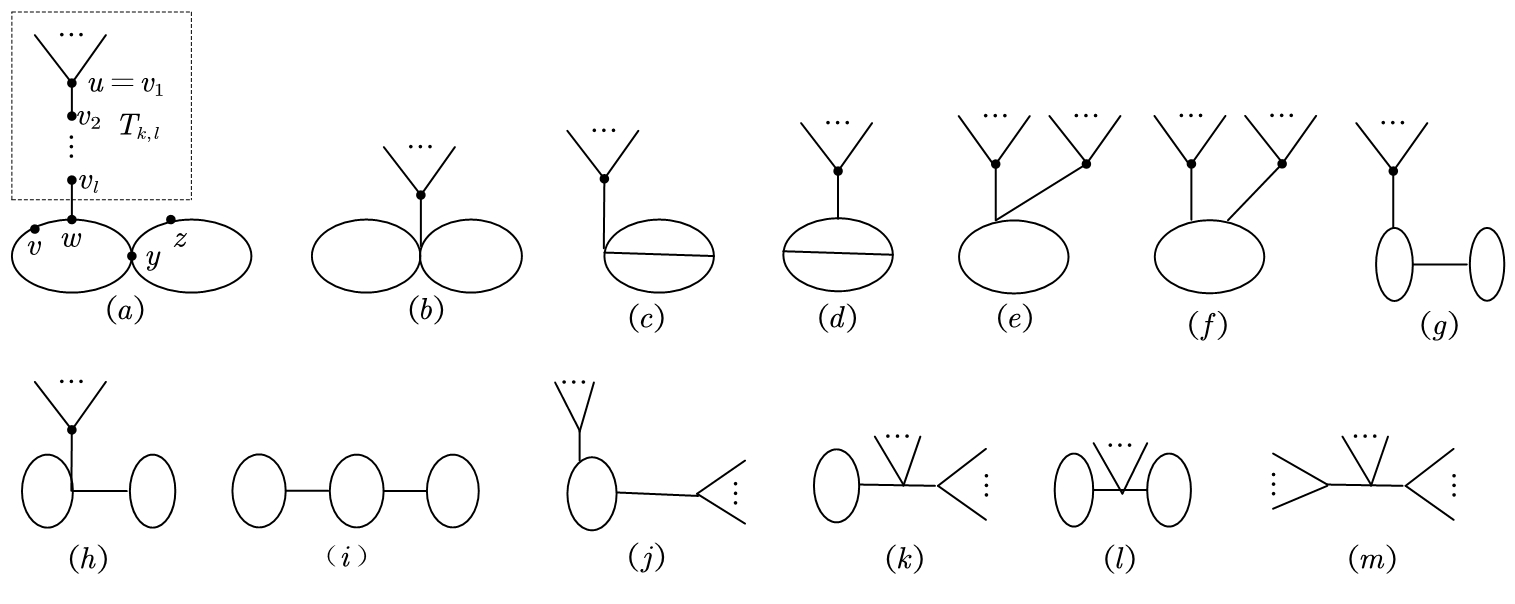}
\caption{$G\in \mathbb{G}_s$}\label{fig1}
\end{figure}

\begin{Lemma}
Let $G\in \mathbb{G}_s$ be one the graph in Fig. 1, then $G$ is not $\lambda$-optimal for $\lambda\notin \sigma(P_s)$.
\end{Lemma}
\begin{proof}
We will prove only case (a), as the proofs for the other cases are similar.
In case (a), as illustrated in Fig. 1, let $u$ be the center vertex of $T_k$, $w$ be a high degree vertex on the cycle $C_1$, $v$ be a neighbor of $w$ on the cycle $C_1$ with $d_G(v)=2$, $y$ be the unique intersection vertex of the two cycles, and $z$ be a neighbor of $y$ on the other cycle $C_2$.
By Lemma 3.1, we have
$$m_G(\lambda) \leq \dim(\mathbb{V}^{\lambda}_G \cup \mathbb{Z}^{\{u,v,z\}}_G) + 3.$$
For $\alpha \in \mathbb{V}^{\lambda}_G \cup \mathbb{Z}^{\{u,v,z\}}_G$, similar to the proof of Lemma 3.7, we can deduce that $\alpha = 0$.
Thus, it follows that $m_G(\lambda) \leq 3 = 2c(G) + q_s(G) - 2$.
Hence, $G$ is not $\lambda$-optimal.
\end{proof}

\begin{Lemma}
Let $G \in \mathbb{G}_s$ and $q_s(G) \geq 1$.
If $\lambda \notin \sigma(P_s)$ and $G$ is $\lambda$-optimal, then $G$ is a graph in which all cycles are disjoint.
\end{Lemma}
\begin{proof}
When $c(G) = 2$ and a series of $s$-$p$-deletion processes is applied to $G$ to obtain $G'\in \mathbb{G}_s$, such that $q_s(G') = 1$, if $G$ is a graph with intersecting cycles, then $G'$ corresponds to one of the cases (a), (b), (c), or (d) in Fig. 1.
Thus, by Lemma 3.9, we conclude that $G'$ is not $\lambda$-optimal.
However, by Lemma 3.7, we also conclude that $G'$ is $\lambda$-optimal, which leads to a contradiction.

When $c(G) \geq 3$, applying a series of $s$-$p$-deletion processes to obtain $G^{\circ}$, if $G$ is a graph with intersecting cycles, according to Theorem 2.1, we know that $G^{\circ}$ is not $\lambda$-optimal by $c(G^{\circ})\geq3$.
This again contradicts Lemma 3.7. Therefore, we conclude that $G$ is a graph in which all cycles are disjoint.
\end{proof}

\begin{Lemma}
Let $G \in \mathbb{G}_s$ with $q_s(G)+c(G)\geq3$, $q_s(G)\geq1$, $c(G)\geq1$ and $\lambda \notin \sigma(P_s)$.
If $G$ is $\lambda$-optimal, then all cycles in $G$ are pendant cycles.
\end{Lemma}
\begin{proof}
Assume that $G$ contains non-pendant cycles.
Since $G$ is a graph in which any two cycles are disjoint, we perform a series of $s$-$p$-deletion and $c$-$p$-deletion processes on $G$ to obtain $G_1\in \mathbb{G}_s$, such that $q_s(G_1)+c(G_1)=3$ and $G_1$ still contains non-pendant cycles.
Combined with Lemma 3.10, $G_1$ can only be one of the configurations (g), (h), (i) or (j) as depicted in Fig. 1.
According to Lemma 3.9, we know that $G_1$ is not $\lambda$-optimal, which contradicts the results established in Lemmas 3.7 and 3.8. Therefore, we conclude that all cycles in $G$ must be pendant cycles.
\end{proof}

\begin{Lemma}
Let $G \in \mathbb{G}_s$ with $q_s(G)+c(G)\geq3$, $q_s(G)\geq1$, $c(G)\geq1$ and $\lambda \notin \sigma(P_s)$.
If $G$ is $\lambda$-optimal, then all $T_{k,l}$ in $G$ are pendant-$T_{k,l}$.
\end{Lemma}
\begin{proof}
Assume that the conclusion is not valid.
We perform a series of $s$-$p$-deletion and $c$-$p$-deletion processes on $G$ to obtain $G_1\in \mathbb{G}_s$, such that $c(G_1) + q_s(G_1)=3$ and $G_1$ contains a non-pendant $T_{k,l}$.
Since every cycle in $G$ is a pendant cycle by Lemma 3.11, $G_1$ can only be one of the configurations (k), (l) or (m) as illustrated in Fig 1.
However, according to Lemma 3.9, we know that $G_1$ is not $\lambda$-optimal, which contradicts the results established in Lemmas 3.7 and 3.8. Thus, the proof is complete.
\end{proof}

\begin{figure}[H]
  \centering
  \includegraphics[width=1.0\linewidth]{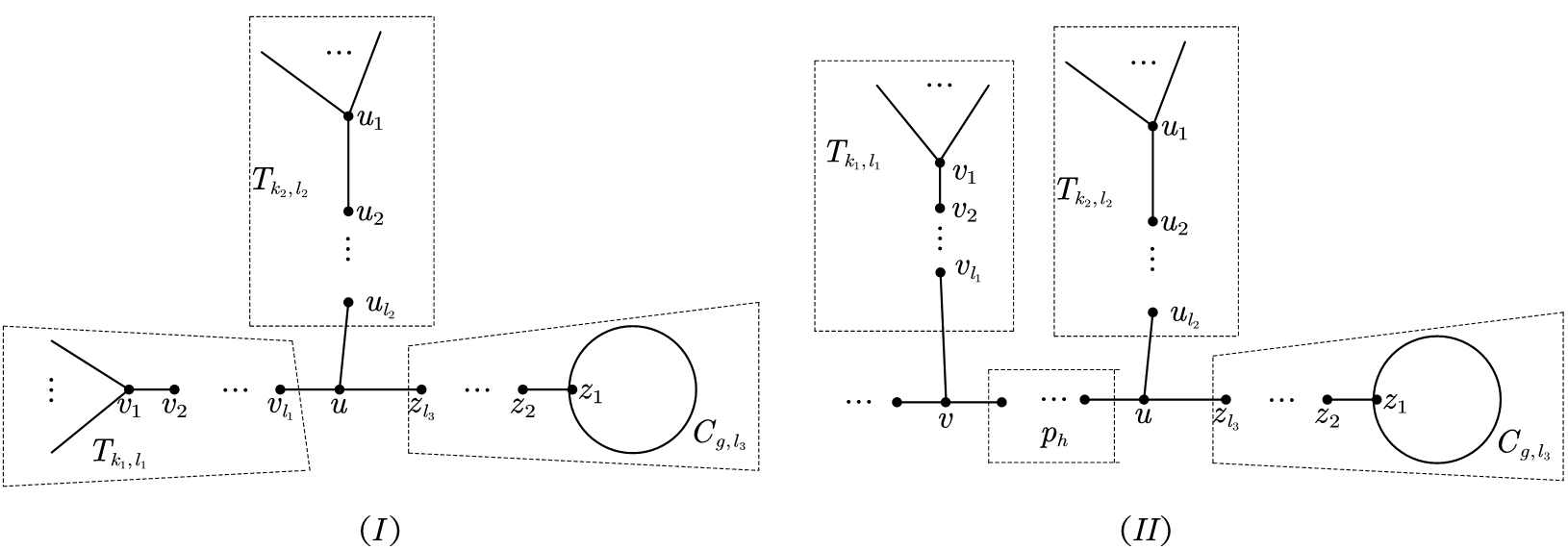}
\caption{$c(G)=1$}\label{fig2}
\end{figure}

\begin{Lemma}
Let $G \in \mathbb{G}_s$ with $q_s(G) \geq 2$, $\lambda \notin \sigma(P_s)$, and $c(G) = 1$.
Then, $G$ is $\lambda$-optimal if and only if $G \in \mathbb{B}(T, 1)$.
\end{Lemma}
\begin{proof}
When $G$ is $\lambda$-optimal, we know from Lemma 3.11 that the only cycle $C$ in $G$ is a pendant cycle.

Next, we will use induction on $q_s(G)$ to prove that $G\in \mathbb{B}(T, 1)$.

When $q_s(G) = 2$,  by combining Lemma 3.11 and 3.12, $G$ is as shown in Fig 2(I), where $z_1 \in V(C) \cap M_G$.
We will prove that $m_{T_{k_1, l_1}}(\lambda)=m_{T_{k_2, l_2}}(\lambda)=1$ and $m_{C_{g,l_3}}(\lambda)=2$.
From Lemma 2.1, we know that
$$m_G(\lambda) \leq \dim(\mathbb{V}^{\lambda}_G \cap \mathbb{Z}^{u}_G) + 1.$$
For $(x_1, x_2, \ldots, x_n)^T = \alpha \in \mathbb{V}^{\lambda}_G \cap \mathbb{Z}^{u}_G$, from $A\alpha = \lambda \alpha$, we can deduce that $\alpha|_{T_{k_i, l_i}} \in \mathbb{V}^{\lambda}_{T_{k_i, l_i}}$ for $1 \leq i \leq 2$.
If $\lambda \notin \sigma(T_{k_1, l_1})$ and $\sigma(T_{k_2, l_2})$, then $\alpha|_{T_{k_i, l_i}} = 0$.
From $A\alpha = \lambda \alpha$ and $\lambda x_i = \sum_{j \sim i} x_j$, we have $\alpha|_{C-z_1} \in \mathbb{V}^{\lambda}_{C-z_1}$.
Thus, $m_G(\lambda) \leq m_{C-z_1}(\lambda) + 1 \leq 2$ by $C-z_1=P_{|C|-1}$, which leads to a contradiction.
Therefore, we can assume  $\lambda$ is an eigenvalue of $T_{k_1, l_1}$.
According to Lemma 3.2, $\lambda$ is not an eigenvalue of $T_{k_1, l_1-1}$.
By Lemma 3.5, we obtain
$$m_G(\lambda)= m_{C_{g, l_3}}(\lambda) + m_{T_{k_2, l_2}}(\lambda) = 2c(G) + q_s(G) - 1 = 3.$$
Furthermore, since $m_{C_{g, l_3}}(\lambda) \leq 2$ and $m_{T_{k_2, l_2}}(\lambda) \leq 1$ by Lemma 3.2 and 3.3, we conclude that $m_{C_{g, l_3}}(\lambda) = 2$ and $m_{T_{k_2, l_2}}(\lambda) = 1$.
Therefore, we have $G \in \mathbb{B}(T, 1)$.

Assuming that the conclusion holds for $q_s(G) = k \geq 2$, we now prove that the conclusion also holds for $q_s(G) = k + 1$.
By performing an $s$-$p$-deletion on $G$, we obtain $G_1=G-T_{k_1,l_1}$ (see Fig. 2(II), let the two pendant-$T_{k_i,l_i}$ that are closest to pendant-$C_{g,l_3}$ be denoted as $T_{k_1,l_1}$ and $T_{k_2,l_2}$).
By Lemma 3.7, we know that $G_1$ is $\lambda$-optimal and that $q_s(G_1) = q_s(G) - 1 = k$.
Therefore, by the inductive hypothesis, we have $G_1 \in \mathbb{B}(T_1,1)$.

Additionally, we similarly perform another $s$-$p$-deletion on $G$ to obtain $G_2=G-T_{k_2,l_2}$.
By the same reasoning, we have $G_2 \in \mathbb{B}(T_2, 1)$.

We will next prove that $G\in \mathbb{B}(T, 1)$.
As illustrated in Fig. 2(II), if $u = v$, then it is evident that $G \in \mathbb{B}(T, 1)$; if $u \neq v$, we only need to demonstrate that the internal path $P_h$ between $u$ and $v$ satisfies $m_{P_h}(\lambda) = 1$.
Since $G_1 \in \mathbb{B}(T_1, 1)$, we can assume $\lambda = 2 \cos \frac{i \pi}{m+1}$ by the definition of $\mathbb{B}(T_1, 1)$, where $i$ and $m + 1$ are coprime, with $1 \leq i \leq m$.
Therefore, we have $l_3 \equiv 1 \ (\bmod \ m+1)$.
Similarly, since $G_2 \in \mathbb{B}(T_2, 2)$, we obtain $l_3 + h + 1 \equiv 1 \ (\bmod \ m+1))$, which implies $h \equiv m \ (\bmod \ m+1))$.
Thus, we conclude that $\lambda \in \sigma(P_h)$.
In summary, this leads us to the conclusion that $G \in \mathbb{B}(T, 1)$.

\end{proof}

\begin{figure}[H]
  \centering
  \includegraphics[width=1.0\linewidth]{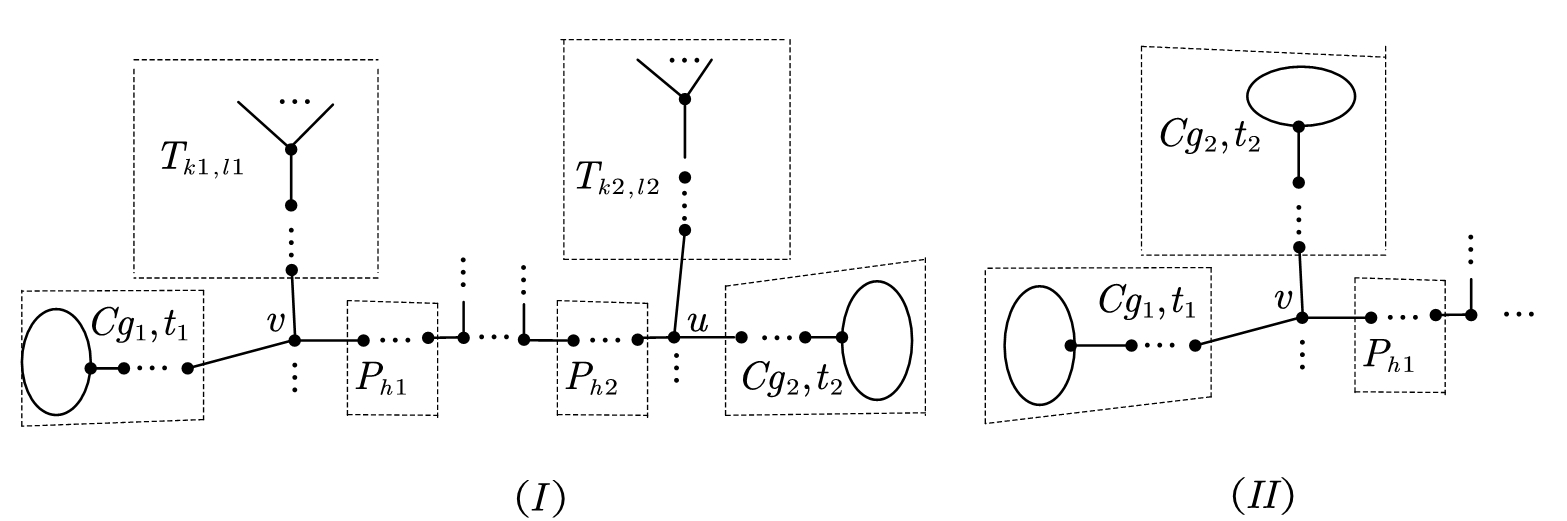}
\caption{$q_s(G)\geq2, c(G)\geq2$}\label{fig3}
\end{figure}

\begin{Theorem}
Let $G \in \mathbb{G}_s$ with $q_s(G) \geq 2$, $\lambda \notin \sigma(P_s)$, and $c(G)\geq1$.
Then, $G$ is $\lambda$-optimal if and only if $G \in \mathbb{B}(T, c(G))$.
\end{Theorem}

\begin{proof}
We perform mathematical induction on $c(G)$.
When $c(G) = 1$, by Lemma 3.13, we know that the conclusion holds.
Assuming that the conclusion is true for $c(G) = k$ with $k \geq 1$, we will prove that it also holds for $c(G) = k + 1$.
We know from Lemma 3.11 and Lemma 3.12 that all cycles and $T_{k_i,l_i}$ in $G$ are pendant.
Therefore, we conduct a $c$-$p$-deletion process on $G$ to obtain $G_1=G-C_{g_1,t_1}$.
According to Lemma 3.8, $G_1$ is $\lambda$-optimal.
Furthermore, since $c(G_1) = c(G) - 1$ and $q_s(G)=q_s(G_1)$, by the inductive hypothesis, we have $G_1 \in \mathbb{B}(T_1, c(G) - 1)$.
On the other hand, we can also perform another $c$-$p$-deletion process on $G$ to obtain $G_2=G-C_{g_2,t_2}$, such that $G_2$ is distinct from $G_1$.
Similarly, we have $G_2 \in \mathbb{B}(T_2, c(G) - 1)$.

Next, we will prove that $G \in \mathbb{B}(T, c(G))$.
We know that the structure of $G$ is represented as in Fig. 3 (I) or (II).
We will only prove the case of (I), while the other case can be proved in a similar manner.
In Fig. 3 (I), if $d_G(u), d_G(v) \geq 4$, then it is evident that $G \in \mathbb{B}(T, c(G))$.
We will discuss three cases:

(1) $d_G(u) = 3$, $d_G(v) \geq 4$;

(2) $d_G(u) \geq 4$, $d_G(v) = 3$;

(3) $d_G(u) = 3$, $d_G(v) = 3$.

We will only prove the case (1), as (2) and (3) can be proved analogously.
When $d_G(u) = 3$ and $d_G(v) \geq 4$, to prove that $G \in \mathbb{B}(T, c(G))$, it suffices to show that $\lambda \in \sigma(P_{h_1})$.
Since $G_1 \in \mathbb{B}(T_1, c(G) - 1)$, we have $m_{T_{k_1, l_1 + 1 + h_1}}(\lambda) = 1$.
Furthermore, since $G_2 \in \mathbb{B}(T_2, c(G) - 1)$, we likewise obtain $m_{T_{k_1, l_1}}(\lambda) = 1$.
By Lemma 3.2, we know that $\lambda \notin \sigma(T_{k_1, l_1 - 1})$.
Thus, by Lemma 3.5, we have
$$m_{T_{k_1, l_1 + 1 + h_1}}(\lambda) = m_{T_{k_1, l_1 - 1}}(\lambda) + m_{P_{h_1}}(\lambda) = 1,$$
which implies $m_{P_{h_1}}(\lambda) = 1$.
Consequently, we conclude that $G \in \mathbb{B}(T, c(G))$, and the proof is complete.
\end{proof}

\begin{figure}[H]
  \centering
  \includegraphics[width=1.0\linewidth]{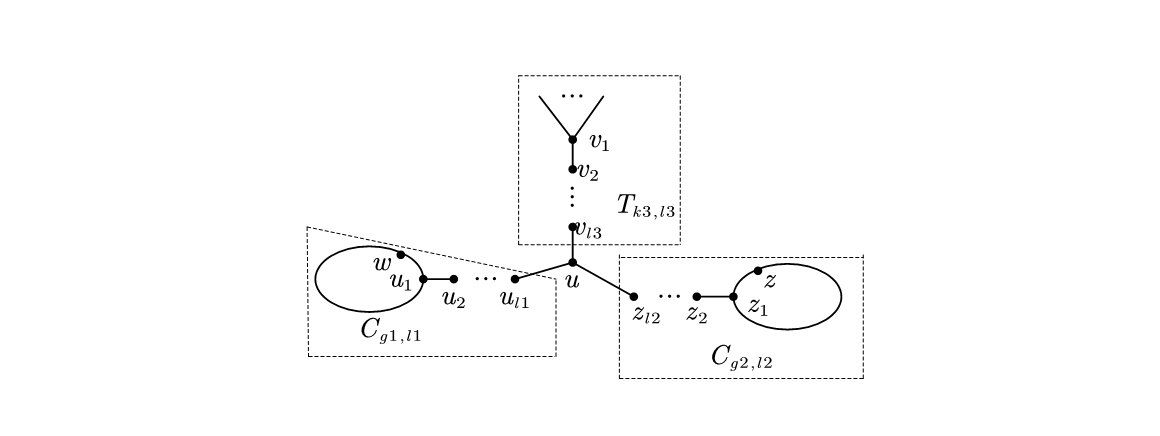}
\caption{$q_s(G)=1$}\label{fig4}
\end{figure}

Next, we will prove the case where $q_s(G) = 1$ and $c(G) \geq 2$.

\begin{Theorem}
Let $G \in \mathbb{G}_s$ with $q_s(G)=1$, $\lambda \notin \sigma(P_s)$, and $c(G)\geq2$.
Then, $G$ is $\lambda$-optimal if and only if $G \in \mathbb{B}(T, q_s(T)-1)=\mathbb{B}(T, c(G))$.
\end{Theorem}

\begin{proof}
We will use mathematical induction to prove our conclusion regarding $c(G)$.
When $c(G) = 2$, $G$ is shaped as shown in Fig. 4.
We aim to prove that $\lambda$ is an eigenvalue of $T_{k_3, l_3}$.
Let $u_1 \in V(C_1) \cap M_G$ and $w$ be a neighbor of $u_1$ on $C_1$.
Similarly, let $z_1 \in V(C_2) \cap M_G$ and $z$ be a neighbor of $z_1$ on $C_2$.
According to Lemma 2.1, we have $m_G(\lambda) \leq \dim(\mathbb{V}^{\lambda}_G \cap \mathbb{Z}^{\{u_1, w, z\}}_G) + 3$.
For $(x_1, x_2, \ldots, x_n)^T = \alpha \in \mathbb{V}^{\lambda}_G \cap \mathbb{Z}^{\{u_1, w, z\}}_G$, from $A\alpha = \lambda \alpha$ and $\lambda x_i = \sum_{j \sim i} x_j$, we can deduce that $\alpha|{T_{k_3, l_3}} \in \mathbb{V}^{\lambda}_{T_{k_3, l_3}}$.
If $\lambda$ is not an eigenvalue of $T_{k_3, l_3}$, then from $\lambda x_i = \sum_{j \sim i} x_j$, we get $\alpha = 0$. Thus, $m_G(\lambda) \leq 3 = 2c(G) + q_s(G) - 2$, leading to a contradiction.
Therefore, $\lambda \in \sigma(T_{k_3, l_3})$.
By Lemma 3.2, we know that $m_{T_{k_3, l_3}}(\lambda) = 1$ and $\lambda \notin \sigma(T_{k_3, l_3-1})$.
By Lemma 3.5, we have
$$m_G(\lambda) = m_{C_{g_1, l_1}}(\lambda) + m_{C_{g_2, l_2}} = 2c(G) + q_s(G) - 1 = 4.$$
Since $m_{C_{g_i, l_i}}(\lambda) \leq 2$ for $1 \leq i \leq 2$, it follows that $m_{C_{g_1, l_1}}(\lambda) = m_{C_{g_2, l_2}}(\lambda) = 2$. Therefore, we conclude that $G \in \mathbb{B}(T, 2) = \mathbb{B}(T, q_s(T) - 1)= \mathbb{B}(T, c(G))$.

Next, we assume that the conclusion holds when $c(G) = k \geq 2$, and we will prove that it also holds when $c(G) = k + 1$.
Since all cycles in $G$ are pendant cycles, we apply a $c$-$p$-deletion process to obtain $G_1$.
According to Lemma 3.8, we know that $G_1$ is $\lambda$-optimal with $c(G_1) = c(G) - 1 = k$.
By the induction hypothesis, we have $G_1 \in \mathbb{B}(T_3, q_s(T_3) - 1)$.

Similarly, we can perform another $c$-$p$-deletion on $G$ to obtain a graph $G_2\neq G_1$, and we can also show that $G_2 \in \mathbb{B}(T_4, q_s(T_4) - 1)$.
Finally, similar to the discussion in Theorem 3.14, we also have $G \in \mathbb{B}(T, c(G))$.
\end{proof}

By combining Theorems 3.14 and 3.15, we can obtain Theorem 2.5.

\section{Acknowledgments}

We gratefully acknowledge the support of the Graduate Innovation Program of China university of Mining and Technology (2024WLKXJ119) and the Postgraduate Research $\&$ Practice Innovation Program of Jiangsu Province (KYCX24{\_}2678) for this work.


\begin{thebibliography}{s2}
\bibitem{CJ1} Q. Chen, J. Guo, Z. Wang, The multiplicity of a Hermitian eigenvalue on graphs. arxiv preprint arxiv:2306.13882, 2023.
\bibitem{LW1} X. Li, Z. Wang, Z. Zhu, Multiplicity of Any Eigenvalue of a Graph with Long Pendant Paths. Available at SSRN 4892226.
\bibitem{RP1} P. Rowlinson, On multiple eigenvalues of trees. Linear Algebra and its Applications, 2010, 432(11): 3007-3011.
\bibitem{SD1} D. Stevanovi$\acute{c}$, Spectral Radius of Graphs, Academic Press, 2015.
\bibitem{TR1} B. Tayfeh-Rezai, The Star Complement Technique. Online Lecture note, URL:http://math. ipm. ac. ir/tayfeh-r/papersandpreprints/starcompl. pdf, 2009.
\bibitem{TK1} K. Toyonaga , C.R. Johnson, Classification of edges in a general graph associated with the change in multiplicity of an eigenvalue. Linear and Multilinear Algebra, 2021, 69(10): 1803-1812.
\bibitem{WL1} L. Wang, L. Wei, Y. Jin, The multiplicity of an arbitrary eigenvalue of a graph in terms of cyclomatic number and number of pendant vertices. Linear Algebra and its Applications, 2020, 584: 257-266.
\bibitem{ZZ1}Y. Zhang, J. Zhao and D. Wong, A characterization for a graph with an eigenvalue of multiplicity $2c(G)+p(G)-1$, Discrete Math. 347 (2024), no.~7, Paper No. 114028, 11 pp.; MR4732854.
\bibitem{ZD1}W. Zhen, D. Wong and Y. Zhang, Eigenvalue multiplicity of graphs with given cyclomatic number and given number of quasi-pendant vertices, Discrete Appl. Math.  347 (2024), 23--29; MR4686339.

\end{thebibliography}
\end{document}